\documentclass[twoside,a4paper]{amsart}
\usepackage{stmaryrd}
\usepackage{ amssymb }
\usepackage{mathabx}
\usepackage{pdfsync}
\usepackage{xcolor}
\definecolor{webblue}{rgb}{0,.5,0}
\definecolor{webblue}{rgb}{.6,0,0}
\definecolor{RoyalBlue}{cmyk}{1, 0.50, 0, 0}
\usepackage[colorlinks=true, breaklinks=true, urlcolor=webblue, linkcolor=RoyalBlue, citecolor=webblue,backref=page]{hyperref}

\newtheorem{theorem}{Theorem}[section]

\newtheorem{corollary}[theorem]{Corollary}

\newcommand{\R}{\mathbb{R}}
\newcommand{\bC}{\mathbb{C}}
\newcommand{\Z}{\mathbb{Z}}
\newcommand{\T}{\mathbb{T}}
\newcommand{\D}{\mathbb{D}}
\newcommand{\N}{\mathbb{N}}

\usepackage{ amssymb }

\title[ Two quantitative versions of Nonlinear Carleson Conjecture]{Two quantitative versions of Nonlinear Carleson Conjecture}
\author{Sergey A. Denisov}

\thanks{
This research was supported by the grants NSF-DMS-2054465, NSF-DMS-2450716, Simons Fellowship in Mathematics, Simons Travel Support for Mathematicians Award, and by the Van Vleck Professorship Research Award. The author gratefully acknowledges the hospitality of IHES where part of this work was done.
}

\address{
\begin{flushleft}
Sergey Denisov: denissov@wisc.edu\\\vspace{0.1cm}
University of Wisconsin--Madison\\  Department of Mathematics\\
480 Lincoln Dr., Madison, WI, 53706,
USA\vspace{0.1cm}\\
\end{flushleft}
}

\begin{document}

\begin{abstract} 
  In this note, we state and compare two quantitative versions of the Nonlinear Carleson Conjecture (NCC). We provide motivations for our conjectures and show that they both imply the NCC. Some applications to the zero distribution of polynomials orthogonal on the unit circle and to their pointwise asymptotics are obtained.
\end{abstract}
\maketitle






\setcounter{tocdepth}{3}

\tableofcontents

We will start with notation: $\N:=\{1,2,\ldots\}, \Z_+:=\{0,1,2,\ldots\}, \T:=\{z\in \bC: |z|=1\}, \D:=\{z\in \bC: |z|<1\}$.
$M(d,\bC)$ stands for the space of the $d\times d$ complex matrices. For $A\in M(d,\bC)$, the symbol $\|A\|_2$ denotes the Frobenius norm: $\|A\|_2:=\sqrt{{\rm tr} (A^*A)}$ and $\|A\|$ denotes the operator norm. For $p\in [1,\infty]$, the symbol $p'$ is the dual exponent: $p'=p/(p-1)$. For a set $S\subset \T$, the symbol $|S|$ indicates its Lebesgue measure.
\section{OPUC, $SU(1,1)$, and matrix products}

Let $\mathcal{M}(\T)$ denote the set of  probabiliy measures  on the unit circle $\mathbb{T}$ whose support is not a finite subset of $\T$. For $\sigma\in \mathcal{M}(\T)$, define
$
F_\sigma(z):=\int_{\mathbb{T}} \frac{\xi+z}{\xi-z}d\sigma,\xi=e^{i\theta},\theta\in [0,2\pi)\,.
$
The function $F_\sigma$ is analytic in $\mathbb{D}$, $\Re F>0$ in $\mathbb{D}$, and $F(0)=1$. The same properties hold for the function $1/F_\sigma$ and therefore \cite{gar} there is a probability measure $\sigma_d\in \mathcal{M}(\T)$ such that
\[
F^{-1}_\sigma(z)=F_{\sigma_d}(z)=\int_{\mathbb{T}}\frac{\xi+z}{\xi-z}d\sigma_d\,.
\]
We will call such a measure $\sigma_d$ dual to $\sigma$. For each $\sigma\in \mathcal{M}(\T)$, one can define the monic orthogonal polynomials $\{\Phi_n(z,\sigma)\}, n\in \Z_+$ on the unit circle (OPUC) as algebraic polynomials that satisfy the following conditions:
\[
\Phi_n(z,\sigma)=z^n+\ldots, \quad \langle \Phi_n,z^j\rangle_\sigma=0,\quad \forall j\in \{0,\ldots,n-1\},\quad \langle f,g\rangle_\sigma:=\int_{\T} f\bar gd\sigma\,.
\]
Denote $\Phi_n:=\Phi_n(z,\sigma), \Psi_n:=\Phi_n(z,\sigma_d)$. For every polynomial $Q$ of degree at most $n$, let $Q^*(z)=z^n\overline{Q(\bar{z}^{-1})}$. Notice that the map $Q\mapsto Q^*$ depends on $n$. It is known \cite{BS} that 
\begin{equation}\label{rec}
\left[\begin{array}{cc}
\Phi_{n+1}&-\Psi_{n+1}\\
\Phi_{n+1}^*&\Psi_{n+1}^*
\end{array}\right]=\left[\begin{array}{cc}
z&-\bar\gamma_n\\
-z\gamma_n&1
\end{array}\right]\left[\begin{array}{cc}
\Phi_{n}&-\Psi_{n}\\
\Phi_{n}^*&\Psi_{n}^*
\end{array}\right], \quad \left[\begin{array}{cc}
\Phi_{0}&-\Psi_{0}\\
\Phi_{0}^*&\Psi_{0}^*
\end{array}\right]=\left[\begin{array}{cc}
1&-1\\
1&1
\end{array}\right]\,,
\end{equation}
where $\{\gamma_g\}$ are the so-called Schur (Verblunsky) parameters and they satisfy $\gamma_g\in \D, \forall g\in \Z_+$. There is a bijection \cite{BS} between the set of measures $\mathcal{M}(\T)$ and the set of Schur parameters $\D^\infty$. If $\{\gamma_g\}$ are Schur parameters for $\sigma$, then the Schur parameters for $\sigma_d$ are $\{-\gamma_g\}$. Denote $\rho_n:=(1-|\gamma_n|^2)^{\frac 12}$ and $\mathcal{R}_n:=\prod_{j\le n}\rho_j, \mathcal{R}_\infty:=\lim_{n\to\infty} \mathcal{R}_n$. Since $\mathcal{R}_n$ is decreasing, $\mathcal{R}_\infty$ always exists and it is zero iff $\{\gamma_g\}\notin \ell^2(\Z_+)$.

The Szeg\H{o} class of measures is defined as the set $Sz(\T):=\{\sigma\in \mathcal{M}(\T): \int_{\T}\log wdm>-\infty\}$ where $m$ is normalized Lebesgue measure on $\T$, i.e., $dm:=d\theta/(2\pi)$, and $d\sigma=wdm+d\sigma_s$ with $\sigma_s$ being the singular measure. The Szeg\H{o} theorem in the theory of OPUC states that $\sigma\in Sz(\T)$ iff $\{\gamma_g\}\in \ell^2(\Z_+)$ and that can be quantified by the following identity \cite{BS}:
\begin{equation}\label{sd2}
 \exp\left(\int_{\T}\log wdm\right)=\prod_{n\ge 0}\rho_n^{2}=\mathcal{R}^2_\infty\,,
\end{equation}
which holds for any $\sigma$ and $\{\gamma_g\}$. Clearly, $\sigma\in Sz(\T)$ iff $\sigma_d\in Sz(\T)$. For measures $\sigma\in Sz(\T)$, the Szeg\H{o} function is defined as $D_\sigma(z):=\exp\left(\frac 12 \int_{\mathbb{T}} \frac{\xi+z}{\xi-z}\log wdm\right)$.
This is the outer function in $H^2(\D)$ that satisfies $D_\sigma(0)>0$ and $|D_\sigma(\xi)|^2=w(\xi)$ for a.e. $\xi\in \T$. 

Sometimes it is more convenient to work with polynomials orthonormal on the unit circle with respect to measure $\sigma$. These polynomials are given by the formula $\phi_n(z,\sigma)=\Phi_n(z,\sigma)/\|\Phi_n(\xi,\sigma)\|_{L^2_\sigma}$. Since (\cite{BS}, formula (1.5.13))
$
\|\Phi_n(\xi,\sigma)\|_{L^2_\sigma}=\mathcal{R}_{n-1}
$,
we get
\[
\sup_{z,n}\left|\frac{\phi_n(z,\sigma)}{\Phi_n(z,\sigma)}-1\right|\lesssim \|\{\gamma_g\}\|_{\ell^2(\Z_+)}^2\,,
\]
when $\|\{\gamma_g\}\|_{\ell^2(\Z_+)}\le \frac 12$. Hence, the problems of studying the size of $\Phi_n$ or $\phi_n$ are identical in the Szeg\H{o} class when  $\|\{\gamma_g\}\|_{\ell^2(\Z_+)}$ is small. For Szeg\H{o} measures, it is known that 
\begin{equation}\label{as}
\lim_{n\to\infty}\Phi^*_n(z,\sigma)=D_\sigma(0)/D_\sigma(z), \quad \lim_{n\to\infty}\phi^*_n(z,\sigma)=1/D_\sigma(z)
\end{equation}
locally uniformly in $\D$ (\cite{BS}, Theorem 2.4.1, p. 144).\medskip

Given any sequence $\{\gamma_n\}: \gamma_n\in \D, n\in \Z_+$, define 
\begin{equation}   \Omega_n:=\rho_n^{-1}\left[\begin{array}{cc}
1&\bar\gamma_n z^{-n}\\
\gamma_nz^n&1
\end{array}\right]\label{matr}
\end{equation}
and $\Pi_n:=\Omega_{n}\cdot \ldots\cdot \Omega_0$. It is a simple fact that $\Pi_n\in SU(1,1)$ for $z\in \T$ and therefore 
\begin{equation}\label{sd8}
\Pi_n=\left[
\begin{array}{cc}
\bar a_n&\bar{b}_n\\
b_n& {a}_n
\end{array}
\right], \quad |a_n|^2=1+|b_n|^2, \quad z\in \T\,.
\end{equation}
The polynomials $a_n$ and $b_n$ have degrees at most $n$, $a_n$ has no roots in $\D$ and $a_n(0)>0$ (see \cite{Kh}, Lemma~4.5 or \cite{nlft}).
If $\gamma_n$ are Schur parameters for $\sigma$, then (see \cite{Kh}, formulas (4.12) and (4.13))
\begin{equation}
\left[\begin{array}{cc}
\Phi_{n+1}&-\Psi_{n+1}\\
\Phi_{n+1}^*&\Psi_{n+1}^*
\end{array}\right]=\mathcal{R}_n \left[
\begin{array}{cc}
z{a}^*_n&-{b}^*_n\\
-zb_n& {a}_n
\end{array}
\right]
\cdot \left[\begin{array}{cc}
1&-1\\
1&1
\end{array}\right]\label{sd4}
\end{equation}
and, therefore,
\begin{equation}\label{sd10}
 a_n=(\Phi^*_{n+1}+\Psi^*_{n+1})/(2\mathcal{R}_n)\,,\quad -zb_n=(\Phi^*_{n+1}-\Psi^*_{n+1})/(2\mathcal{R}_n)\,.
\end{equation}
Denote
\begin{equation}\label{dfg1}
o_1(\xi):=\sup_{n\ge 0}|a_n(\xi)-1|,\quad o_2(\xi):=\sup_{n\ge 0}|b_n(\xi)|, \quad O(\xi):=\sup_{n\ge 0}\left\|\Pi_n(\xi)-I\right\|_2, \quad \xi\in \T
\end{equation}
and we clearly have 
\begin{equation}\label{sd1}
\|\Pi_n-I\|_2^2=2(|a_n-1|^2+|b_n|^2)\le 2( o_1^2+o_2^2).
\end{equation}\medskip

\section{Two versions of quantified NCC and motivations for them}

The famous Carleson-Hunt theorem \cite{Swed} says that 
\begin{equation}\label{ch}
\|M_0f\|_{L^p[0,1)}\le_{p}\|f\|_{L^p[0,1)}, \quad 
M_0f:=\sup_{N\in \mathbb{N}}\left|\sum_{|n|<N}\widehat f_ne^{2\pi i nx}\right|
\end{equation}
for $p\in (1,\infty)$, where $\widehat f_n=\int_0^1 fe^{-2\pi i nx}dx$ is the Fourier coefficient of a function $f$. Function $M_0$ is an associated maximal function. The bound \eqref{ch} was used to settle Lusin's conjecture, i.e., it showed that {\it for every $f\in L^p(0,1), p\in (1,\infty]$, we have $\lim_{N\to\infty}\sum_{|n|<N}\widehat f_ne^{2\pi i nx}=f(x)$ for a.e. $x\in [0,1)$.}
The analog of Lusin's conjecture in the OPUC theorem is called Nonlinear Carleson Conjecture (NCC) and it says: {\it is it true that asymptotics \eqref{as} holds for a.e. $z\in \T$ assuming that $\sigma\in Sz(\T)$?} This problem has a long history and it can also be formulated for Dirac equations, Schr\"odinger operators and Krein systems (nonlinear Fourier transform)  \cite{ck,tao2,tao3,pks,nlft}. The general theory of the nonlinear Fourier transform is an active field with multiple applications \cite{michel1,liban,kovac,sylv,nlft}. One purpose of the current note is to put forward two quantifications of such a conjecture,  motivate them and compare. For $\sigma\in \mathcal{M}(\T)$, define the maximal function  $M(\xi,\sigma):=\sup_{n}|\Phi^*_n(\xi,\sigma)-1|, \xi\in \T$. From \eqref{sd10}, we get  
\begin{equation}\label{sd11}
o_{1(2)}(\xi)\le (M(\xi,\sigma)+M(\xi,\sigma_d)/(2\mathcal{R}_\infty)+\mathcal{R}_\infty^{-1}-1, \quad \xi\in \T\,.
\end{equation}

\noindent {\bf Conjecture, [qNCC-I].} There is $\epsilon>0$ such that $\|\{\gamma_g\}\|_{\ell^2(\Z_+)}\le \epsilon$ implies
\begin{equation}
\int_{\T} M^2(\xi,\sigma)d\sigma\lesssim \|\{\gamma_g\}\|^2_{\ell^2(\Z_+)}\,.
\end{equation}

\noindent {\bf Conjecture, [qNCC-II].} There is $\epsilon>0$ such that $\|\{\gamma_g\}\|_{\ell^2(\Z_+)}\le \epsilon$ implies
\begin{equation}\label{vtor}
\int_{\T} \log (1+O^2(\xi,\{\gamma_g\}))dm\lesssim \|\{\gamma_g\}\|^2_{\ell^2(\Z_+)}\,.
\end{equation}

 We will see later that the assumption on $\ell^2$-norm of $\{\gamma_g\}$ being small is not restrictive when studying the problem of pointwise asymptotics.
Our first result compares these two conjectures.

\begin{theorem}
We have qNCC-I $\Rightarrow$ qNCC-II.
\end{theorem}

\begin{proof}Assume qNCC-I holds. Due to \eqref{sd1}, we have $\log (1+O^2)\lesssim \log(1+o_1^2+o_2^2)$ and it is sufficient to show that
\begin{equation}\label{sd16}
\int_{\T} \log (1+o_1^2+o_2^2)dm\lesssim \|\{\gamma_g\}\|^2_{\ell^2(\Z_+)}\,.
\end{equation}
Now, \eqref{sd11} implies $\log (1+o_j^2(\xi)) \lesssim  \log(1+M^2(\xi,\sigma))+\log(1+M^2(\xi,\sigma_d))+\|\{\gamma_g\}\|_{\ell^2(\Z_+)}^2, j\in \{1,2\}$.
We can write
\begin{align}\label{sd100}
\int_{w>1}\log(1+M^2(\xi,\sigma))dm\le \int_{w>1}\log(1+M^2(\xi,\sigma)w)dm\le\int \log(1+M^2(\xi,\sigma)w)dm\\\le
 \int M^2(\xi,\sigma)wdm\stackrel{{\rm (qNCC-I)}}{\lesssim} \|\{\gamma_g\}\|^2_{\ell^2(\Z_+)}\,.\label{sd13}
\end{align}
Then,
\begin{align*}
\int_{w\le 1}\log(1+M^2(\xi,\sigma))dm=
-\int_{w\le 1} \log wdm+\int_{w\le 1}\log wdm+\int_{w\le 1}\log(1+M^2(\xi,\sigma))dm=\\
-\int_{w\le 1}\log wdm+\int_{w\le 1}\log(w+M^2(\xi,\sigma)w)dm\le\\ -\int_{w\le 1}\log wdm+\int_{w\le 1}\log(1+M^2(\xi,\sigma)w)dm
\stackrel{\eqref{sd2}+\eqref{sd100}+\eqref{sd13}}{\lesssim}  \|\{\gamma_g\}\|_{\ell^2(\Z_+)}^2\,.
\end{align*}
By qNCC-I applied to $\sigma_d$, we get
$\int_{\T} M^2(\xi,\sigma_d)w_ddm\lesssim  \|\{\gamma_g\}\|_{\ell^2(\Z_+)}^2
$ and we similarly have \[\int  \log(1+M^2(\xi,\sigma_d))dm\lesssim \|\{\gamma_g\}\|^2_{\ell^2(\Z_+)}\]  for $\|\{\gamma_g\}\|_{\ell^2(\Z_+)}$ small enough. Hence, \eqref{sd16} follows.
\end{proof}\medskip
Next, we discuss the motivations for our conjectures. The  Menshov-Rademacher theorem  \cite{Kash,MT} states:
\begin{theorem}[Menshov-Rademacher] Suppose $\{\chi_n(\xi)\},n\ge 1$ is an orthonormal system in $L^2_\sigma(\T)$, then
\[
\int_{\T} \sup_{n}\left|\sum_{j=1}^n \alpha_j\chi_j(\xi)\right|^2d\sigma\lesssim \sum_{j\ge 1}|\alpha_j\log (1+j)|^2
\]
for every sequence $\{\alpha_j\}$.
\end{theorem}
The connection between qNCC-I and the previous result is almost immediate (see \cite{pks}, Section 8) if we take the recursion $\Phi_{n+1}^*=\Phi_n^*-\gamma_nz\Phi_n$ and $\Phi^*_0=1$ into account (see \eqref{rec}). Indeed, 
\[
\Phi_n^*(z,\sigma)=1-z\sum_{j=0}^{n-1}\gamma_j\|\Phi_j\|_{2,\sigma}\phi_j(z,\sigma)
\]
so
\[
M(\xi,\sigma)\le\left|\sup_{n\ge 1}\sum_{j=0}^{n-1}\gamma_j\|\Phi_j\|_{2,\sigma}\phi_j(\xi,\sigma)\right|
\]
and $ \|M\|_{2,\sigma}\lesssim\Bigl(\sum_{j\ge 0}|\gamma_j|^2\log^2 (2+j)\Bigr)^{\frac 12}$. The qNCC-I suggests that for OPUC the logarithm in the last sum can be dropped provided that $\|\{\gamma_g\}\|_{\ell^2(\Z_+)}$ is small.\medskip

We continue with the motivation for qNCC-II.
Take $d\ge 2$ and consider $\gamma(t): [0,1] \mapsto GL(d,\bC)$, a smooth curve. Let $\gamma_r(t):=\int_0^t \gamma'\gamma^{-1}d\tau:   [0,1] \mapsto M(d,\bC) $ be its right trace. One can define the distance between $A,B\in GL(d,\bC)$  by
\[
d(A,B):=\inf_{\gamma(0)=A, \gamma(1)=B}\int_0^1 \|\gamma'\gamma^{-1}\|dt\,.
\] 
Given $q\in [1,\infty]$, the variational norm of a continuous curve $\Gamma: [0,1]\mapsto M(d,\bC)$ is defined as
\[
\|\Gamma\|_{V^q}:=\left\{
\begin{array}{cc}
\sup_{n\in \N}\sup_{0=t_0<t_1<\ldots<t_n=1}\left(\sum_{j=0}^{n-1}\|\Gamma(t_{j+1})-\Gamma(t_j)\|^q\right)^{1/q}, &q<\infty,\\
{\rm diam}\,(\Gamma), &q=\infty
\end{array}
\right.\,.
\]
In \cite{tao3}, Lemma C.3., the authors proved, in particular, that
\begin{equation}\label{sd15}
\|\gamma\|_{V^q}\le \|\gamma_r\|_{V^q}+C_q\min (\|\gamma_r\|_{V^q}^q, \|\gamma_r\|_{V^q}^2)\quad \|\gamma_r\|_{V^q}\le \|\gamma\|_{V^q}+
C_q\min(\|\gamma\|_{V^q}^q, \|\gamma\|_{V^q}^2)
\end{equation}
 for every $q\in (1,2)$ which gives
\begin{equation}\label{sd19}
\|\gamma\|_{V^q}\le \|\gamma_r\|_{V^q}+C_q\|\gamma_r\|_{V^q}^2, \quad \|\gamma_r\|_{V^q}\le \|\gamma\|_{V^q}+
C_q\|\gamma\|_{V^q}^2
\end{equation}
for all $q\in (1,2)$.

  In \cite{tao3}, the authors use \eqref{sd15} along with a variational version of Menshov-Paley-Zygmund theorem to improve earlier results of Christ and Kiselev (see \cite{ck,cknew} and \cite{ck8}, formula (1.3)).
Now, consider the sequence $\{\Upsilon_n\}, \Upsilon_n\in GL(d,\bC), n\in \Z_+$ where 
$
\sup_{n}\|\Upsilon_n-I\|\le\frac 12
$. One can define the piecewise smooth curve $\gamma(t), t\ge 0$ as the solution to $\gamma'=V\gamma, \gamma(0)=I$ where
\[
V(t):=\log \Upsilon_j, \quad t\in [j,j+1), \, j\in \Z_+\,.
\]
We notice that $V(t)=\Upsilon_j-I+\Delta_j, t\in [j,j+1)$ where $\|\Delta_j\|\sim \|\Upsilon_j-I\|^2$.
Then, $\Upsilon_n\cdot\ldots\cdot \Upsilon_0=\gamma(n+1)$ and \eqref{sd15} can be applied to the product of matrices (the assumption made in \cite{tao3}, Lemma C.3 that the curve $\gamma$ is smooth can be relaxed to piecewise smooth). 
Recall that 
$
{\rm dist} (I,A)\sim \log(1+\|I-A\|)
$ for $A\in SU(1,1)$.
Taking $\Upsilon_n=\Omega_n(\xi,\{\gamma_g\})$ as in \eqref{matr}, applying \eqref{sd19} to $\gamma$, and the variational Menshov-Paley-Zygmund theorem for Fourier series (see section B in \cite{tao3} for Fourier integral version) gives us, in particular (compare with p. 461, \cite{tao3}),
\[
\|\log^{\frac 12}(1+O^2)\|_{L^{p'}(\T)}\le_p \|\{\gamma_g\}\|_{\ell^p(\Z_+)}
\]
for $p\in (1,2)$.
Setting $p=2$ in this estimate gives the conjectured bound \eqref{vtor}.\medskip

\section{Some applications of qNCC}

We have
\begin{theorem}\label{t67} Suppose $\{\gamma_g\}\in \ell^2(\Z_+)$ and qNCC-II holds, then 
$\lim_{n\to\infty}\Pi_n(\xi)$ exists for a.e. $\xi\in \T$ and 
$\lim_{n\to\infty}\phi_n^*(\xi,\sigma)=D_\sigma^{-1}(\xi)$ for a.e. $\xi\in \T$.
\end{theorem}
\begin{proof}Assume  qNCC-II holds. We will show that the sequence $\{\Pi_n(\xi,\{\gamma_g\})\}$ is Cauchy for a.e. $\xi\in \T$.  We have
$
\|\Pi_{n+p}-\Pi_n\|=\|(\Omega_{n+p}\cdot \ldots \cdot \Omega_{n+1}-I)\Pi_n\|\le O\cdot  \|(\Omega_{n+p}\cdot \ldots \cdot \Omega_{n+1}-I)\|\,.
$
Since $O(\xi)<\infty$ for a.e. $\xi\in \T$,
we only need to show that 
\[
\limsup_{N\to\infty} \,O_N(\xi)=0,\quad O_N:=\sup_{n\ge N,p\ge 1}\|\Omega_{n+p}\cdot \ldots \cdot \Omega_{n+1}-I\|
\]
for a.e. $\xi\in \T$. Since \[
\int \log (1+O_N^2)dm\stackrel{({\rm qNCC-II})}{\lesssim} \|\{\gamma_g\}\|^2_{\ell^2(\ge N)},
\]
we can apply Markov's inequality:
for every $\epsilon>0$ and $N\in \N$, we have
\[
|\{\xi\in \T: O_N(\xi)>\epsilon\}|\lesssim \frac{\|\{\gamma_g\}\|^2_{\ell^2(\ge N)}}{\log(1+\epsilon^2)}\,.
\]
Since $\{O_N\}\searrow \limsup_{N\to\infty}O_N$ as $N\to\infty$, we have $|\{\xi\in \T: \limsup_{N\to\infty}O_N(\xi)\ge \epsilon\}|=0$ for every $\epsilon>0$. Therefore, $\limsup_{N\to\infty} O_N(\xi)=0$ for a.e. $\xi\in \T$ and we have convergence of $\lim_{n\to\infty}\Pi_n$ by applying the Cauchy criterion. The formula \eqref{sd4} provides convergence of $\{\phi^*_n(\sigma,\xi)\}$ for a.e. $\xi\in \T$. Since (see \cite{Kh}, Corollary 5.11)
\[
\lim_{n\to\infty}\int \left|\frac{1}{\phi_n^*}-D_\sigma\right|^2dm=0\,,
\]
we get $\lim_{n\to\infty}\phi_n^*(\xi,\sigma)=D_\sigma^{-1}(\xi)$ a.e.
\end{proof}

It is known that all zeroes of $\phi_n(z,\sigma)$ are inside $\D$. The following theorem shows, in particular, that for the Szeg\H{o} measures, the pointwise convergence of $|\phi_n|$ on $\T$ is equivalent to the condition that its zeroes stay away from $\T$. To state this result, we need some notation first. 
Given a parameter $\rho\in (0,1)$ and a point $\xi\in \T$, define the Stolz angle $S_\rho^*(\xi)$ to be the convex hull of $\rho\D$ and $\xi$. Let $\{f_n\}$ denote the Schur interates (see \cite{Kh}) for the measure $\sigma$.  
\begin{theorem}[Bessonov-Denisov, \cite{BD}] \label{t2}
Let $\sigma \in Sz(\T)$ and  $Z(\phi_n) = \{z \in \D: \; \phi_n(z,\sigma) = 0\}$. Take any $a>0$ and denote $r_{a,n}= 1-a/n$. Then, for almost every $\xi\in \T$, the following assertions are equivalent:
\begin{itemize}
\item[$(a)$] $\lim_{n \to \infty} |\phi_{n}^{*}(\xi)|^2= |D_{\sigma}^{-1}(\xi)|^2$,
\item[$(b)$] $\lim_{n \to \infty} {\rm dist}(Z(\phi_n), \xi)\,n = +\infty$,
\item[$(c)$] $\lim_{n \to \infty} f_n(r_{a,n}\xi) = 0$,
\item[$(d)$] $\lim_{n\to \infty}\sup_{z\in S_\rho^*(\xi)}|f_n(z)|=0$ for every $\rho\in (0,1)$.
\end{itemize} 
\end{theorem}
Our previous results imply the following corollary.
\begin{corollary}Suppose that either $\{\gamma_g\}\in \ell^p(\Z_+)$ for some $p\in [1,2)$ or $\{\gamma_g\log(g+2)\}\in \ell^2(\Z_+)$, then (a)-(d) from the previous theorem hold for a.e. $\xi\in \T$.
\end{corollary}\medskip

\section{Connections to Carleson-Hunt maximal function and $SU(1,1)$ version of a theorem of Calderon and Stein}
In conclusion, we put forward the weakened versions of qNCC-I and qNCC-II:
\[
\int \sup_{n\ge 0}|\Phi_n(\xi,\sigma)|^2d\sigma\le 1+\nu_1(\|\{\gamma_g\}\|_{\ell^2(\Z_+)}) \quad {\rm  (wqNCC-I)}
\]
and
\[
\int  \, \log\, \sup_{n\ge 0}\|\Pi_n(\xi,\{\gamma_g\})\|dm\le \nu_2(\|\{\gamma_g\}\|_{\ell^2(\Z_+)})\quad {\rm  (wqNCC-II)}
\]
assuming that $\|\{\gamma_g\}\|_{\ell^2(\Z_+)}<\epsilon$ with some $\epsilon>0$ and $\nu_1$, $\nu_2$ are certain functions that satisfy $\lim_{t \downarrow 0}\nu_{1(2)}(t)=0$. It is not hard to see that qNCC-I implies wqNCC-I,  qNCC-II implies wqNCC-II, and wqNCC-I implies wqNCC-II. The converse statements are not known. It is also not known if these weak versions imply the NCC. The inequality in \eqref{sd1} shows that $O^2\ge 2o_2^2$ and, therefore,
$\int \log(1+2o_2^2)dm\le \int \log(1+O^2)dm$. The next result shows that the weakened version of qNCC implies the known Carleson-Hunt bound \eqref{ch}  for $p=2$. 
\begin{theorem}If there is $\epsilon>0$ so that
\begin{equation}\label{loi1}
\int \log(1+o_2^2)dm\lesssim \|\{\gamma_g\}\|^2_{\ell^2(\Z_+)}
\end{equation}
for every sequence $\{\gamma_g\}: \|\{\gamma_g\}\|_{\ell^2(\Z_+)}\le \epsilon$, 
then 
\begin{equation}\label{pol3}
\int \sup_{n}\Bigl|\sum_{j\le n}\gamma_je^{ij\theta}\Bigr|^2dm\lesssim  \|\{\gamma_g\}\|^2_{\ell^2(\Z_+)}
\end{equation}
 for every $\{\gamma_g\}\in\ell^2(\Z_+)$.
\end{theorem}
\begin{proof}Fix any $\{\gamma_g\}\in \ell^2(\Z_+)$. The map $b_n(\xi,\{\lambda\gamma_j\}): (\xi,\lambda,\{\gamma_j\}_{ j\le n})\in \T\times\D\times  \D^{n+1}\mapsto \bC $ satisfies \cite{nlft} the bound
\[
\left|b_n(\xi,\{\lambda\gamma_j\})-\lambda\sum_{j\le n}\gamma_j\xi^j\right|\le C(\{\gamma_g\})|\lambda|^2
\]
as $\lambda\to 0$. Given any $N\in \N$ and a measurable map $N(\xi): \xi\in \T\mapsto \{0,\ldots,N\}$, we get 
\[
\int \log (1+|b_{N(\xi)}(\xi,\{\lambda\gamma_j\})|^2)dm\lesssim |\lambda|^2\|\{\gamma_g\}\|^2_{\ell^2(\Z_+)}
\]
from \eqref{loi1}. As $\lambda\to 0$, we get
\[
|\lambda|^2\int \Bigl|\sum_{j\le N(\xi)}\gamma_j\xi^j\Bigr|^2dm+O(|\lambda|^3)\lesssim |\lambda|^2\|\{\gamma_g\}\|^2_{\ell^2(\Z_+)}\,.
\]
Hence,
\[
\int\Bigl|\sum_{j\le N(\xi)}\gamma_j\xi^j\Bigr|^2dm\lesssim \|\{\gamma_g\}\|^2_{\ell^2(\Z_+)}\,.
\]
Since $N$ and $N(\xi)$ are arbitrary, we obtain \eqref{pol3}.\end{proof}

The following result is classical (see \cite{guzman}, Section 2.1) and has many generalizations \cite{stein}. It predates Carleson's proof of Luzin's conjecture and is attributed to Calderon and to Stein.

\begin{theorem}\label{carl}The following statements are equivalent:
\begin{itemize}

\item[A.] For every $f\in L^2(0,1)$, the sequence $\{\sum_{j=-n}^n \widehat f_je^{2\pi ijx}\}$ of partial Fourier sums converges for a.e. $x\in [0,1)$

\item[B.] We have the weak  $(2,2)$ type estimate for the maximal function $M_0$, i.e.,
\[
|\{x\in [0,1): |(M_0f)(x)|\ge \lambda\}|\lesssim \frac{\|f\|_2^2}{\lambda^2}
\]
for all $\lambda>0$ and all $f\in L^2(0,1)$.
\end{itemize}
\end{theorem}
Below, we will adjust its proof to obtain an $SU(1,1)$ nonlinear version.
We list some properties of the transform $\{\gamma_g\}\mapsto \{\Pi_n(\xi,\{\gamma_g\})\}$  we need. They are immediate from the definition \eqref{sd8}. Firstly, for every $\eta\in \T$, we have
\begin{equation}\label{smur1}
\Pi_n(\xi,\{\eta\gamma_g\})=\left[
\begin{array}{cc}
\bar a_n&\overline{\eta{b}_n}\\
\eta b_n& {a}_n
\end{array}
\right], \quad \|\Pi_n(\xi,\{\eta\gamma_g\})-I\|_2=\|\Pi_n(\xi,\{\gamma_g\})-I\|_2\,.
\end{equation}
Secondly, suppose $\xi_0=e^{i\alpha_0}$. Then,
\begin{equation}\label{shift}\text{\it rotation by $\alpha_0$}:\quad 
\Pi_n(\xi,\{\gamma_g\xi_0^{-g}\})=\Pi_n(\xi/\xi_0,\{\gamma_g\}).
\end{equation}
For each $N\in \N$, define the $\{\gamma^{[N]}_g\}$ by 
\[
\gamma^{[N]}_g:=\left\{
\begin{array}{cc}
0, &g<N,\\
\gamma_{g-N}, &g\ge N
\end{array}\right..
\] The last property we need is  
\[
\text{\it the right shift by $N$ coordinates:}\quad 
\Pi_{n+N}(\xi,\{\gamma^{[N]}_g\})=\left[
\begin{array}{cc}
\bar a_n&\bar{b}_n \xi^{-N}\\
b_n\xi^N& {a}_n
\end{array}
\right]\,.
\]
Notice that 
\begin{equation}
\left\|\left[\begin{array}{cc}
\bar a_n&\bar{b}_n \xi^{-N}\\
b_n\xi^N& {a}_n
\end{array}\right]-I\right\|_2=\left\|\left[\begin{array}{cc}
\bar a_n&\bar{b}_n \\
b_n& {a}_n
\end{array}\right]-I\right\|_2
\end{equation}
and, therefore, for each $\eta,\xi_0\in \T$, we get
\begin{equation}\label{spas}
\|\Pi_{n+N}(\xi,\{\eta\gamma^{[N]}_g\xi_0^{-g}\})-I\|_2\stackrel{\eqref{smur1}+\eqref{shift}}{=}\|\Pi_{n}(\xi/\xi_0,\{\gamma_g\})-I\|_2\,.
\end{equation}
Consider the following function (recall \eqref{dfg1} and our second conjecture in \eqref{vtor})
\[
G(\alpha,\beta):=\sup_{\{\gamma_g\}:\|\{\gamma_g\}\|_{\ell^2(\Z_+)}\le \beta} |\{\xi\in \T:  \log(1+O^2(\xi,\{\gamma_g\}))\ge \alpha\}|\,.
\]
Clearly, $G\le |\T|=2\pi$, it is  increasing in $\beta$ and  decreasing in $\alpha$. 
\begin{theorem}
Assume that for every $\{\gamma_g\}: \|\{\gamma_g\}\|_{\ell^2(\Z_+)}\le\frac 12$, the sequence $\{\Pi_n(\xi,\{\gamma_g\})\}$ converges for  $\xi\in W(\{\gamma_g\})\subset \T$, where the set $W(\{\gamma_g\})$ can depend on $\{\gamma_g\}$ and has positive Lebesgue measure. Then, 
\begin{equation}\label{kot1}
G(\alpha,\beta)\le_{\alpha}\beta^{2}
\end{equation}
for every  $\beta\le \frac 12$.  
Conversely, if \eqref{kot1} holds for every $\alpha>0$ and $\beta\le \frac 12$, then the sequence $\{\Pi_n(\xi,\{\gamma_g\})\}$ converges for a.e. $\xi\in \T$.
\end{theorem}
\begin{proof}
Recall that the matrices $\Pi_n$ in \eqref{sd8} satisfy $\|\Pi_n\|=\|\Pi_n^{-1}\|$ since $\Pi_n\in SU(1,1)$.   Therefore, we have
\begin{equation}\label{liu3}
\lim_{N\to\infty}\sup_{n\ge N, p\in \N}\|\Omega_{n+p}(\xi,\{\gamma_g\})\cdot\ldots\cdot \Omega_{n+1}(\xi,\{\gamma_g\})-I\|=0
\end{equation}
for each $\xi\in W(\{\gamma_g\})$. We continue the proof  by assuming that \eqref{kot1} fails for some $\alpha$. Then, choosing the subsequence if necessary, we can find a sequence $\{\beta_j\}\downarrow 0, \beta_j\le 2^{-j}$ such that 
\begin{equation}\label{kot2}
G(\alpha,\beta_j)\ge \beta_j^{2}j^2\,.
\end{equation}
Hence, for each $j$, there is $\{\gamma^{(j)}_g\}$ with finite support, i.e., 
$\gamma^{(j)}_g=0$ for all $g\ge N_j$, such that $\|\gamma^{(j)}_g\|_{\ell^2(\Z_+)}\le \beta_j$ and 
\begin{equation}\label{leb2}
|E_j|\ge \beta_j^{2}j^2/2\,,
\end{equation}
 where
\begin{equation}\label{lash1}
E_j:=\{\xi\in \T:  \log(1+O^2(\xi,\{\gamma^{(j)}_g\}))\ge \alpha/2\}\,.
\end{equation}
Notice that, since $\gamma^{(j)}_g=0$ for $g\ge N_j$, we get  $\sup_{n\ge 0} \log(1+\|\Pi_n(\xi,\{\gamma^{(j)}_g\})-I\|_2^2)=\sup_{n< N_j}  \log(1+\|\Pi_n(\xi,\{\gamma^{(j)}_g\})-I\|_2^2)$. We claim that there is a sequence $\{d_j\}$ of natural numbers such that 
\begin{equation}\label{liu4}
\sum_{j\ge 1} d_j\beta_j^2<\infty, \quad \sum_{j\ge 1} d_j\beta_j^{2}j^2=\infty\,.
\end{equation}
Indeed, it suffices to chose $d_j$ such that $d_j\sim \beta_j^{-2}j^{-2}$ and recall that $\beta_j\le 2^{-j}$. We consider the sequence of sets $\{E_p^*\}$ in which each set $E_j$ is repeated $d_j$ times, i.e.,
\[
\{E_p^*\}:=\{
\underbrace{E_1, \ldots, E_1}_{d_1}, \underbrace{E_2,\ldots, E_2}_{d_2}, \ldots\}\,,
\]
so that $\sum_{p\ge 1} |E^*_p|\stackrel{\eqref{leb2}+\eqref{liu4}}{=}\infty$ and, therefore (see Lemma 2.1.2 in \cite{guzman}), there is a sequence  $\{\xi^*_p\},\,\, \xi^*_p=e^{i\alpha_p^*}\in \T$ such that the shifted sets $\widehat E_p:=\{e^{i(\alpha_p^*+\beta)}:  e^{i\beta}\in E^*_p\}$ satisfy $|\limsup_p \widehat E_p|=|\T|=2\pi$, where $\limsup_p\widehat E_p:=\cap_{j\in \N}\cup_{s\ge j}\widehat E_s$.

Our next goal is to use $\{\gamma^{(m)}_g\}, m\in \N$ to construct $\{\gamma^*_g\}\in \ell^2(\Z_+)$ for which the assumption of the Theorem is violated.
 Let $d_0:=0, N_0:=0, P_{j}:=N_0d_0+\ldots+N_jd_j, j\ge 0$. Then, for each $\widehat E_p$ with $p=d_0+\ldots+d_j+s, 1\le s\le d_{j+1}$, we let 
 \[
 \gamma^*_{\ell+(s-1)N_{j+1}+P_{j}}:=\gamma^{(j+1)}_\ell\cdot \xi_p^{-\ell},  \quad 0\le \ell\le N_{j+1}-1\,.
 \]
 For such a choice, we can use  \eqref{spas} and \eqref{lash1} to obtain
 \[
 \sup_{0\le \ell\le N_{j+1}-1} \|\Omega_{\ell+(s-1)N_{j+1}+P_{j}}(\xi,\{\gamma^*_g\})\cdot\ldots\cdot \Omega_{(s-1)N_{j+1}+P_{j}}(\xi,\{\gamma^*_g\})-I\|_2\ge (e^{\alpha/2}-1)^{\frac 12}
 \]
for $\xi\in \widehat{E}_p$ and every $p$. Since $\|\{\gamma^*_g\}\|_{\ell^2(\Z_+)}\stackrel{\eqref{liu4}}{<}\infty$, we have a contradiction with \eqref{liu3} for  $\xi\in W(\{\gamma^*_g\})\cap \limsup_p\widehat E_p\neq \emptyset$. The first claim of the theorem is proved. The proof of the second one is standard and repeats the argument from the proof of Theorem \ref{t67}.
\end{proof}



\bibliographystyle{plain}

\bibliography{samplebib}

\end{document}